\documentclass[a4paper, fleqn, 12pt]{article}
\usepackage[a4paper]{geometry}
\usepackage[utf8]{inputenc}
\usepackage[english]{babel}
\usepackage[T1]{fontenc}
	\usepackage{multirow}
\usepackage{ytableau}
\usepackage{booktabs}
\usepackage{makecell}
\usepackage{mathtools}
\usepackage{amssymb}
\usepackage{paralist}
\usepackage{tikz}
\usepackage{amsmath}
\usepackage{mathrsfs}
\usepackage{tikz-cd}
\usepackage{amsthm}
\usepackage{lscape}
\usepackage{enumitem}
\usetikzlibrary{arrows}

\newtheorem{theorem}{Theorem}[section]
\newtheorem{lemma}[theorem]{Lemma}
\newtheorem{corollary}[theorem]{Corollary}

\newtheorem{proposition}[theorem]{Proposition}

\theoremstyle{definition}

\newtheorem{remark}[theorem]{Remark}
\newtheorem{definition}[theorem]{Definition}

\numberwithin{subcase}{case}
\DeclareMathOperator{\Chi}{{\mathfrak X}}
\DeclareMathOperator{\diag}{diag}

\DeclareMathOperator{\ind}{ind}
\DeclareMathOperator{\Irr}{Irr}
\DeclareMathOperator{\id}{id}
\DeclareMathOperator{\GL}{GL}

\DeclareMathOperator{\Gal}{Gal}
\DeclareMathOperator{\Aut}{Aut}
\DeclareMathOperator{\End}{End}
\DeclareMathOperator{\disc}{disc}
\DeclareMathOperator{\Fix}{Fix}
\DeclareMathOperator{\Br}{Br}

\DeclareMathOperator{\GO}{GO}
\DeclareMathOperator{\GU}{GU}
\DeclareMathOperator{\SO}{SO}
\DeclareMathOperator{\SU}{SU}
\DeclareMathOperator{\SL}{SL}

\newcommand{\disj}{\stackrel{.}{\cup}}
\newcommand{\Disj}{\stackrel{.}{\bigcup}}
\newcommand{\g}{{\bf g}}
\newcommand{\h}{{\bf h}}
\newcommand{\z}{{\bf z}}
\newcommand{\w}{{\bf w}}
\newcommand{\T}{{\bf t}}
\newcommand{\Z}{\mathbb{Z}}
\newcommand{\Q}{\mathbb{Q}}
\newcommand{\Res}{\mathrm{Res}}
\newcommand{\F}{\mathbb{F}}

\begin{document}
\title{Unitary Discriminants of $\SL_3$(q) and $\SU_3$(q)}
\author{Linda Hoyer\footnote{linda.hoyer@rwth-aachen.de} \ and Gabriele Nebe\footnote{nebe@math.rwth-aachen.de}}
\date{Lehrstuhl f\"ur Algebra und Zahlentheorie, RWTH Aachen University, Germany}
\maketitle
\begin{abstract}
We give a full list of the unitary discriminants of the even degree indicator 'o' ordinary irreducible characters of $\mathrm{SL}_3(q)$ and $\mathrm{SU}_3(q)$.
\end{abstract}

\section{Introduction}

Let $G$ be a finite group, $K$  a number field, $n\in \mathbb{N}$ and let $\rho : G\to \GL_n(K)$ be a $K$-representation of $G$ of degree $n$. Then $\rho$ is uniquely determined by its character $\chi : G\to K, g \mapsto \mathrm{trace}(\rho(g))$. The ordinary character table for $G$ lists the values of 
all irreducible characters on the conjugacy classes of $G$. 
Together with an additional number-theoretic invariant, the Brauer element of $\chi $ (see Definition \ref{Brauer}), 
it contains all necessary information about the linear actions of $G$ over 
number fields, i.e. the representations $\rho $ as above.

As $G$ is finite, the image $\rho (G) \subseteq \GL_n(K)$  is contained in either a symplectic, unitary, or orthogonal group. 
The papers \cite{Nebe2022OrthogonalS} and \cite{Unitary} develop 
methods towards classifying the orthogonal and unitary groups that contain $\rho(G)$.
These methods are then applied
to compute the orthogonal discriminants of 
the even degree indicator $+$ irreducible characters 
and the unitary discriminants of the even degree indicator 'o' characters
of a large portion of the groups in the ATLAS of finite groups \cite{ATLAS}.

The ATLAS of finite groups contains the character tables of 
small finite simple groups, including all sporadic simple groups. 
Most of the non-abelian simple groups are finite groups of Lie Type.
They fall in infinite series for which there are so called generic character tables
parameterising the representations. 
For these infinite series, it is necessary to compute the corresponding generic orthogonal and unitary discriminants.

The first formula in the literature for such generic discriminants is 
the one by Jantzen and Schaper (see for instance \cite{JamesMathas}),
giving the composition factors
of the discriminant groups of the Specht modules for all symmetric groups.
The first generic orthogonal character tables have been obtained in 
\cite{SL2} for the groups $\SL_2(q)$ and all prime powers $q$. 
Later the authors of this paper computed the generic orthogonal discriminants 
of the groups $\SL_3(q)$ and $\SU_3(q)$ \cite{SL3SU3}, again for all 
prime powers $q$. 

This paper continues this investigation by determining the 
unitary discriminants of the groups $\SL_3(q)$ and $\SU_3(q)$ for all prime powers $q$.  
The computation illustrates three different methods: 
For $\SL_3(q)$ Harish-Chandra induction from the parabolic subgroup 
$\GL_2(q) \ltimes \mathbb{F}_q^2$ is enough to obtain all unitary discriminants
(Theorem \ref{SL3A} and \ref{SL3B}). 

For the unitary groups we also 
distinguish between even and odd defining characteristic. 
Section \ref{SU3gen} lists the important players and the relevant facts 
that apply to both situations. 
For $2$-powers $q$ the irreducible even degree indicator 
'o' characters of $\SU_3(q)$ have degree $q(q^2-q+1)$. 
These characters appear in a rank 2 monomial representation 
that can be analysed using similar methods as for rank 2 permutation representations
(see Section \ref{SU3even}). 
The last section deals with unitary group in odd defining characteristic. 
Here we need to use the full strength of the methods from \cite{Unitary} 
and a metabelian subgroup of the Borel subgroup to finally deduce the 
results for odd $q$. 

This paper is the starting point of a long term project to 
determine the unitary discriminants of ordinary irreducible characters of finite 
groups of Lie Type. On the one hand, the explicit results obtained in this paper
allow to obtain unitary discriminants for characters that are Harish-Chandra 
induced from parabolic subgroups with Levi factor of Type $U_3$ or $L_3$. 
On the other hand, the methods applied here can be generalised to obtain 
(perhaps less explicit) results for finite groups of Lie Type of higher rank. 

This paper is supported by the Deutsche Forschungsgemeinschaft
(DFG, German Research Foundation),
Project-ID 286237555.

\section{Methods}

\subsection{Quadratic and Hermitian forms} 
Let $L$ be a field of characteristic $\neq 2$ 
and $\sigma \in \Aut(L)$ an automorphism of order $1$ or $2$. 
Put $K:=\Fix_{\sigma}(L) $ to denote the fixed field of $\sigma $. 
An $L/K$ Hermitian space is a finite dimensional $L$-vector space $V$ 
together with a non-degenerate $L/K$-Hermitian form 
$H:V\times V \to L$, i.e. a $K$-bilinear map
such that $H(av,w) = aH(v,w)$ and $H(v,w)= \sigma(H(w,v))$ for
all $v,w\in V$ and $a\in L$.
Let
$N:= N_{L/K}(L^{\times })$ denote the norm subgroup of $K^{\times }$.
Then $$(K^{\times })^2 \leq N \leq K^{\times } $$
and $N=(K^{\times })^2$ if $\sigma = \id $.
In the latter case, $L=K$ and $H$ is usually called a symmetric bilinear 
form. 

\begin{definition}
        The {\em discriminant} of $H$ is the class
	$$ \disc(H) :=   (-1)^{{{n}\choose{2}}}\det (H_B) N \in K^{\times } / N,$$
	where $n=\dim (V)$ and 
        $H_B := (H(b_i,b_j))_{i,j =1 }^n \in L^{n\times n}$ is the Gram matrix of
        $H$ with respect to any basis $B=(b_1,\ldots , b_n)$ of $V$.
\end{definition}

For $L\neq K$ the class of the discriminant 
defines a unique quaternion algebra over $K$ which
we call the discriminant algebra of $H$:

 \begin{definition}
	 \begin{itemize}
		 \item[a)] 
For $a,b\in K^{\times }$ let $(a,b)_{K} $ denote the central simple
        $K$-algebra with $K$-basis $(1,i,j,ij)$ such that
        $i^2=a$, $j^2=b$, $ij=-ji$.
		 \item[b)] 
        For a quadratic extension $L= K[\sqrt{\delta }]$ of $K$ we put
        $$(L,b)_K  := (\delta,b)_{K} .$$
	Then the class of $[(L,b)_{K}] \in \Br(L,K) $ lies 
         in the 
         Brauer group  of central simple $K$-algebras
        that are split by $L$.
	We denote by 
			 $$\disc_L([(L,b)_K]) := b N_{L/K}(L^{\times})$$
			 the {\em $L$-discriminant} of $[(L,b)_K]\in \Br(L,K)$.
\item[c)]
      For $\sigma \neq \id$   the {\em discriminant algebra} of the Hermitian
         form $H$ with discriminant $dN$ 
        is defined as the class
         $\Delta (H) := [(L,d)_{K}] \in \Br(L,K) .$
	 \end{itemize}
\end{definition}

Note that $\Delta (H) $ is the Clifford invariant of 
the quadratic form $Q_H:V_K \to K , v \mapsto H(v,v) $
(see for instance \cite[Remark (10.1.4)]{Scharlau}). 
This result also implies the well definedness of the 
discriminant algebra, which 
 is also guaranteed by 
the following remark.

\begin{remark}  \label{quatsplit}
For a quadratic extension $L/K$ and $a,b\in K^{\times }$ we have 
$(L,a)_K \cong (L,b)_K$ if and only if $aN_{L/K}(L^{\times }) = b N_{L/K}(L^{\times })$.
\end{remark}

\begin{remark}\label{Ldiscdiv}
Division algebras over number fields $K$ have the local-global property: 
They are isomorphic if and only if they become isomorphic over all 
	completions of $K$. In particular, we can identify a quaternion algebra ${\mathcal Q}$
by listing all the places $\wp_1,\ldots , \wp_h$ 
	of $K$, where the completion is still a division algebra.
	The primes $\wp_1,\ldots, \wp_h$ are the {\em ramified primes} in 
	${\mathcal Q}$ and we put
	$${\mathcal Q}=: {\mathcal Q}_K( \wp_1,\ldots , \wp_h ) .$$
\end{remark}

The finite ramified places in ${\mathcal Q}$ are exactly those that divide the 
discriminant ideal of a maximal $\Z_K$-order in ${\mathcal Q}$. 
	If $b\in {\mathbb{Z}} _K\setminus \{0\} $  then
$\Z_L \oplus \Z_L j $ is a $\Z_K$-order in $(L,b)_K$ of discriminant 
$\disc_K(L)^2 b^2 $. In this situation we hence have the following remark:

\begin{remark}\label{Ldiscdiv2}
	Let $K$ be a number field and assume that $b\in {\mathbb{Z}} _K\setminus \{0\} $. 
	Then any  be a finite place $\wp $ that ramifies in $(L,b)_K$ 
	divides 
	$\disc_K(L) b .$
\end{remark}

It is a general principle that odd degree extensions are not relevant for 
isometry of quadratic or Hermitian forms. 

\begin{lemma}\label{odd} (see \cite[Corollary 6.16]{BOI})
Let $F$ be an odd degree extension of $K$. Put $E:=FL$ and 
	extend $\sigma $ to a field automorphism $\sigma $ of $E$ 
	with fixed field $F$. 
	Let $(V,H)$ and $(W,H')$ be two $L/K$ Hermitian spaces.
	If the $E/F$ Hermitian spaces $(V\otimes_L E,H_E) $ and 
	$(W\otimes _L E, H'_E)$ are isometric, then also 
	$(V,H) \cong (W,H')$.
\end{lemma}

\subsection{Unitary stable characters}

For a finite group $G$ let $\Irr(G)$ denote the set of  irreducible 
complex characters of $G$. The {\em Frobenius-Schur indicator} $\ind(\chi )$
of $\chi \in \Irr(G)$ takes values in $\{ + , -, \mbox{'o'} \}$. 
We have $\ind(\chi ) = $ 'o' if and only if the character field ${\mathbb Q}(\chi )$ 
is not real, $\ind(\chi ) = +$ if there is a real representation affording 
the character $\chi $, and $\ind(\chi ) = -$ for real characters $\chi $ that are 
not afforded by a real representation.
We denote by 
$$\Irr^+(G):= \{ \chi \in \Irr(G) \mid \ind (\chi ) = +, \ \chi(1) \mbox{ even } \} $$ 
the even degree indicator $+$ characters of $G$ and by 
$$\Irr^o(G):= \{ \chi \in \Irr(G) \mid \ind (\chi ) = \mbox{'o'} , \  \chi(1) \mbox{ even } \} $$ 
the even degree indicator 'o' characters of $G$.

\begin{definition}(\cite[Definition 5.12]{Nebe2022OrthogonalS})
	A character $\chi $ of a finite group $G$ is called 
	{\em orthogonal} if there is a representation $\rho $ affording the
	character $\chi $ and fixing a non-degenerate quadratic form. 
	An orthogonal character is called {\em orthogonally stable}
	if and only if all its indicator $+$ constituents have even degree.
\end{definition}

Then \cite[Theorem 5.15]{Nebe2022OrthogonalS} 
shows that the orthogonally stable characters 
are exactly those orthogonal characters that have a well defined discriminant: 

\begin{definition}
	Let $\chi $ be an orthogonally stable character of a finite group $G$ 
	with character field $K:={\mathbb Q}(\chi )$.
	Then the {\em orthogonal discriminant} $\disc(\chi ) \in K^{\times}/(K^{\times})^2$ is the unique square class of the character field such that for any orthogonal representation $\rho : G\to \GL_n(L)$ and any non-degenerate $\rho(G)$-invariant quadratic form 
	$Q$ we have that $\disc(Q) = \disc(\chi ) (L^{\times })^2 $.
\end{definition}

The paper \cite{Unitary} transfers the notion of orthogonal stability to 
the Hermitian case: 

\begin{definition}\label{unitarystable} \cite[Definition 5.10]{Unitary}
        An ordinary character $\chi $ of a finite group $G$
        is called {\em unitary stable} if
        all irreducible  
        constituents of $\chi $ have even degree.
\end{definition}

Note that a unitary stable orthogonal character is orthogonally stable,
 but the converse is not always the case.
Similarly as in the orthogonal case we get 
that a unitary stable character has a well defined unitary discriminant: 

\begin{remark} (\cite[Proposition 5.13]{Unitary})
	Let $\chi $ be a unitary stable character of the finite group $G$ 
	and assume that the 
	character field $L:={\mathbb Q}(\chi )$ is not real. Denote by $K$  
 its real subfield. 
	Then there is a unique class 
	$$[\Delta (\chi )] \in \Br(L,K) $$ 
	such that for all representations $\rho :G \to \GL_n(M)$ and all 
	$\rho (G)$-invariant non-degenerate Hermitian forms $H$ 
	we have $$\Delta(H) = [\Delta (\chi) \otimes_K M^+] $$
	where $M^+$ is the maximal real subfield of $M$.
	\\
	$[\Delta (\chi )] $ is called the {\em unitary discriminant algebra} of $\chi$.
\end{remark}

Whereas orthogonal discriminants of orthogonally stable characters 
are essentially independent of the chosen splitting field, 
for unitary discriminants this field matters. 

\begin{definition}
	Let $\chi $ be a unitary stable character of the finite group $G$ 
	and let $L$ be a totally 
	complex number field with real subfield $K\neq L$.
	Assume that there is an $L$-representation $\rho $ of 
	$G$ affording the character $\chi $. 
	Then all $\rho(G)$-invariant non-degenerate Hermitian forms $H$
	have the same discriminant $\disc(H) =: d N_{L/K}(L^{\times })$. 
	Then $\disc_L(\chi) := d N_{L/K}(L^{\times })$ is called the 
	{\em $L$-discriminant} of $\chi  $ and 
	$\Delta _L(\chi ):= (L,d)_{K}$ the 
	{\em $L$-discriminant algebra} of $\chi $.
\end{definition}

\subsection{Induced representations} 

Many irreducible characters $\chi $ of finite groups $G$ are imprimitive, i.e. induced from 
a character $\psi $ of a proper subgroup $U$. 
Then a $G$-invariant form in the induced representation 
is just the orthogonal sum of $[G:U]$ copies of 
a $U$-invariant form.
However, the character field of $\psi $ might be larger than the one of $\chi $,
and we only get the discriminant over the character field of $\psi $ 
(see \cite[Remark 7.2]{Unitary}). 
In view of Lemma \ref{odd} it is helpful to know when 
$[{\mathbb Q}(\psi ): {\mathbb Q}(\chi )]$ is odd:

\begin{lemma} \label{induz}
	Let $G$ be a finite group, $U\leq G$, $\psi \in \Irr(U)$ such that 
	$\chi:=\psi^G \in \Irr(G)$. 
	Then $\mathbb{Q}(\chi ) \leq \mathbb{Q}(\psi )$. 
	Let $\Psi := \{  \psi _1 ,\ldots , \psi _h \} $ be the 
	set of constituents of 
 the restriction $\chi _{|U }$ of $\chi $ to $U$ 
	of  degree $\psi_i(1) = \psi (1)$ and assume that the cardinality, $h$, 
	of $\Psi $ is odd. Then there is $i_0 \in \{ 1,\ldots , h \}$ such that 
	 $[\mathbb{Q}(\psi _{i_0}) : \mathbb{Q}(\chi )] $ is odd.
\end{lemma} 

\begin{proof}
	By Frobenius reciprocity  $\Psi = \{ \psi _i \in \Irr (U) \mid \chi = \psi _i^G \}$. 
For $\psi _i \in \Psi $ a full regular orbit under the Galois group 
	$\Gal (\mathbb{Q}(\psi _{i}) / \mathbb{Q}(\chi )) $ 
	is contained in $\Psi $ and $\Psi $ is a disjoint union of 
	such Galois orbits.
	As $|\Psi |$ is odd at least one of these orbits has odd length.
\end{proof}

\subsection{Schur indices} 

Schur indices play an important role in the computation of unitary discriminants 
of unitary stable characters. As described in Theorem \ref{fixalg} below we need all 
local Schur indices to compute these discriminants.
These are encoded in the Brauer element, a notion coined in 
\cite[Definition 2.1]{Turull}:

\begin{definition} \label{Brauer}
	Let $\chi $ be an irreducible ordinary character 
	of some finite group $G$. 
	Let $K$ be some field containing the 
	character field ${\mathbb{Q}}(\chi )$  and 
         let $\rho $ be a 
	$K$-representation of $G$ affording the character  $m\chi $ for some positive integer $m$.
	Then $m$ is minimal if and only if $\End_{KG}(\rho )  =: D$ is 
	a central simple division algebra over $K$. 
	In this case $m^2 = \dim_K(D)$ and $m_K(\chi):=m $ 
	is the {\em Schur index} of $\chi $ over $K$.
	The class $[\chi ]_K:=[D] \in \Br(K)$ of $D$ in the Brauer group $\Br(K)$ of $K$ 
	is called the {\em Brauer element} of $\chi $ over $K$.
	\\
	If $K=\Q(\chi )$ is the character field of $\chi$, 
	then we sometimes omit the field, so $[\chi ] := [\chi ]_{\Q(\chi )}$.
	\\
	The field $K$ is called a {\em splitting field} of $\chi $ 
	if the Brauer element $[\chi ]_K = [K]$ is  trivial.
\end{definition}

To compute unitary discriminants we sometimes need to compute local Schur indices of 
certain characters. General computational methods are described in 
\cite{Unger} and \cite{Lorenz}. When other 
methods fail they all fall back to the following 
very useful result from \cite{Benard}: 

\begin{theorem}\label{Benard} \cite[Theorem 8.1]{Benard} 
	Let $\chi $ be an ordinary irreducible character of a finite group 
	lying in a $p$-block with cyclic defect group. 
	Let $\varphi $ be a $p$-modular constituent of $\chi $. 
	Then the $p$-adic Schur index of $\chi $ is 
	$$m_{\mathbb{Q}_p(\chi )}(\chi ) = [ \mathbb{Q}_p (\chi ,\varphi ) : \mathbb{Q}_p (\chi )] .$$
\end{theorem}

Explicit computations of Schur indices for monomial characters can be obtained from 
\cite[Theorem 2]{Yamada} which we recall for the reader's convenience.

\begin{theorem} \cite[Theorem 2]{Yamada}  \label{Yamada}
	Let $K$ be a field of characteristic $0$.
	Let $H$ be a normal subgroup of a finite group $G$ and let 
	$\psi $ be a linear character of $H$ such that the induced character 
	$\chi:= \psi ^G$ is irreducible. 
	Then the character field $K(\chi )$ is a subfield of $K(\psi )$.
	Put $\Gamma := \Gal (K(\psi ),K(\chi))$ and let 
	$$F:= \{ \g \in G \mid \psi^\g  = \psi ^{\gamma(\g ) } \mbox{ for some } \gamma(\g ) \in \Gamma \} .$$ 
	Then $\gamma :F/H \cong \Gamma $ is a group isomorphism. 
	Choose coset representatives $\g_1,\ldots , \g_n$ of $F$ in $H$,
	put $\gamma_i := \gamma(\g_iH)$ and 
	let $\h_{ij}\in H$ such that $\g_i \g_j = \g_k \h_{ij} $ for some $k\in \{1,\ldots ,n\}$. Then 
	$$\beta : \Gamma \times \Gamma \to K(\psi ) , 
	\beta (\gamma _i, \gamma_j ):= \psi(\h_{ij}) $$ 
	is a factor system and $[\chi ]_{K(\chi)} $ is 
	the inverse of the class of the crossed product algebra defined by $\beta $.
\end{theorem}

\subsection{Unitary discriminants of real characters} \label{indicator}

Let $\chi $ be an irreducible real character of $G$ of even degree 
$2m=\chi(1)$. 
Then its Frobenius-Schur indicator is either $-$ or $+$ and the Brauer element 
$[\chi ]$ is represented by  a quaternion algebra 
${\mathcal Q}$ over the totally real number field 
$K=\Q(\chi )$.
The following results hold in more generality
 (see \cite{Unitary}), but we only use them here for splitting fields
 $L$ of $\chi $ so for simplicity we 
assume that $L $ is a totally complex quadratic extension 
of $K$ that is a maximal subfield of ${\mathcal Q}$.

\begin{proposition} \label{ind+-} 
	\begin{itemize} 
		\item[a)] (\cite[Proposition 5.12]{Unitary})
	If the indicator of $\chi $ is $-$, then 
			$$\disc _L(\chi) = \disc_L([\chi]_K)^{m} \mbox{ and } 
			\Delta (\chi ) = [\chi]_K^m .$$
			If $m$ is even then $\disc_L(\chi )$ is 
			trivial.
		\item[b)] (\cite[Proposition 5.12]{Unitary}, \cite[Proposition 10.33]{BOI})
			Assume that  the indicator of $\chi $ is $+$. Then the following hold:
			\\
			(i) The character 
		$\chi $ is orthogonally stable and  has 
		a well defined orthogonal discriminant $\disc(\chi )$. 
		\\
			(ii) The orthogonal discriminant $\disc (\chi )$ is represented
			by $(-1)^m$ times the reduced norm of a skew symmetric unit 
			in any representation affording $\chi $ 
			(see \cite[Proposition 2.2]{NebeOrtDet}).
			\\
			(iii) 
			The $L$-discriminant 
			and the discriminant algebra of $\chi $ are obtained as
		$$\disc _L(\chi) = \disc(\chi ) \disc_L([\chi ]_K)^{m} \mbox{ and } 
		\Delta (\chi ) = [(L, \disc(\chi ))_K ][\chi]_K^m .$$
			(iv) 
			If $m$ is even or $[\chi]_K = [K]$, 
			then $\disc_L(\chi ) = \disc(\chi )$.
	\end{itemize} 
\end{proposition}

Combining  Proposition \ref{ind+-} with \cite[Theorem 4.7]{Albanian} 
yields an easy formula for the unitary discriminant of a unitary stable rational
character of a $2$-group.

\begin{corollary} \label{2groups} 
	Let $G$ be a $2$-group and $\chi \in \Irr(G)$ be an irreducible 
	rational character of degree $\chi (1)$ a multiple of $4$. 
	Then the unitary (or orthogonal) discriminant of $\chi $ is $1$ over
	any splitting field of $\chi $.
\end{corollary}

\subsection{Orthogonal subalgebras}  \label{OS} 

The paper \cite{Unitary} transfers the computational methods from 
\cite{Nebe2022OrthogonalS} to the Hermitian case. 
One important notion is the one of orthogonal subalgebras. 
For  $\chi \in \Irr^o(G)$ 
 the following method is quite useful: 
Put $L:={\mathbb Q}(\chi )$ and assume that there is an $L$-representation
$\rho : G \to \GL_{2m}(L) $
affording the character $\chi $.
Let $K$ denote the real subfield of $L$ and let
$H$ be a non-degenerate $\rho(G)$-invariant $L/K$-Hermitian form.
Assume that there is $\alpha \in \Aut (G)$
with $\alpha ^2 = 1$ such that $\chi \circ \alpha  = \overline{\chi }$.
Then $\chi + \overline{\chi }$ extends to an 
irreducible character
$\Chi = {\textrm{Ind}}_G^{\mathcal G}(\chi ) $
 of the semidirect product 
${\mathcal G}:= G \rtimes \langle \alpha \rangle $.

Let
$A:=\Fix_{\alpha }(\rho ) = \langle \rho (g) + \rho(\alpha (g)) \mid g\in G \rangle _K $
denote the $\alpha $-fixed algebra.
Then $A$ is invariant under the adjoint involution  of $H$.
Denote the restriction of this involution to $A$ by $\iota _A$.

\begin{theorem} (\cite[Theorem 9.1]{Unitary})  \label{fixalg}
\begin{itemize}
\item[(a)] 
        $[A] = [\Chi ]_K \in \Br(L,K)$.
        \item[(b)]
Assume that
 the Frobenius-Schur indicator of $\Chi $ is $+$.
 \\
		(i) 
                Then $\iota _A$ is an orthogonal involution on $A$.
		Its  discriminant
$\disc(\iota _A)$ is $(-1)^m$ times the
square class in $K$ of the reduced norm of any skew symmetric unit $X$ in $A$
(see \cite[Proposition 2.2]{NebeOrtDet}). 
\\
		(ii) 
 The discriminant of $H$ is $\disc(H) = \disc_L([A])^m \disc(\iota _{A} )$.
        \item[(c)]
 Assume that the Frobenius-Schur indicator of $\Chi $ is $-$. \\ Then
		$\disc(H) = \disc_L([\Chi]_K)^m$ and
                $\Delta (H) = [\Chi]_K^m$.
        \end{itemize}
\end{theorem}

\subsection{A generalisation of the $Q_8$-trick} 

Sometimes subgroups help us to predict all unitary discriminants of 
faithful characters. The following result is a generalisation of 
the $Q_8$-trick \cite[Section 8.1]{Unitary}. 
A variant is later used in Theorem \ref{GL2}.

\begin{theorem}\label{extraspecial}
	Let $d\geq 2$ and $m:=2^d a$ for some $a\in {\mathbb{N}}$.
	Let $p,\ell $ be not necessarily distinct odd primes and 
	let $q:=p^f$ be some power of $p$. 
	Let $G$ be some subgroup of $\GL_m(q)$ containing 
$\SO_m^+(p)$ 
	(e.g. $\SL_m(q) \leq G \leq \GL_m(q), 
	 \SU_m(q) \leq G \leq \GU_m(q), \SO_m^+(q) \leq G \leq \GO_m^+(q) $).
	 Let $\chi $ be an irreducible faithful 
	 ordinary or $\ell $-Brauer character of $G$. 
	 Then
	 \begin{itemize}
	 \item[(a)] The character degree $\chi(1)$ is a multiple of $2^d$. 
	 \item[(b)] If the indicator of $\chi $ is $+$, then its orthogonal 
		 discriminant is $1$.
	 \item[(c)] If $\chi $ is an ordinary character of indicator 'o', then
		 its unitary discriminant is $1$. 
	 \end{itemize}
\end{theorem}

\begin{proof}
	Consider the group $U:=2^{1+2d}_+ \cong \otimes ^d D_8$. 
	Then $U$ has a unique irreducible character
	$\psi $ that restricts non-trivially to the center $Z(U) \cong C_2$.
	This character has degree $\psi (1) = 2^d$, indicator $+$, 
	and is the character of an integral irreducible 
	representation $\rho: U\to \SL_{2^d}(\mathbb{Z}) $.  
	Moreover, $\psi $ is orthogonally stable of orthogonal discriminant 
	$1$. 
	Reducing $\rho $ mod $p$ hence shows that $U
	\leq \SO_m^+(p) \leq \SL_m(p) $, so $U\leq G$. 
	As $\chi $ is a faithful character of $G$, its restriction 
	to $U$ is a multiple of $\psi $. In particular, $2^d = \psi (1)$ 
	divides $\chi (1)$. Also (b) and (c) follow from the 
	fact that $\chi _{|U } = a \psi $ is orthogonal and unitary stable. 
\end{proof}

\section{The special linear group}

Let $p$ be a prime and let $q$ be a power of $p$. 
The group $\mathrm{SL}_3(q)$ is the group of $3 \times 3 $-matrices over 
$\F _q$ of determinant 1. 
It contains a maximal parabolic subgroup 
$$P:=\left\{ \left( \begin{array}{ccc} a & b & c \\ 0 & d & e \\
0 & f & g \end{array} \right) \in \SL_3(q) \right\} $$
of odd index $[\SL_3(q) : P ] = q^2+q+1 $. 
Note that  $P$ is the semidirect product $P \cong \GL_2(q) \ltimes \mathbb{F}_q^2  $
where the action of $ \GL_2(q)$ on $\mathbb{F}_q^2$ is given by
$$\g \cdot \h = \det(\g) (\g \h) \mbox{ for } \g \in \GL_2(q), \h \in \mathbb{F}_q^2 .$$
The center of $P$ is the center $Z$ of $\SL_3(q)$, hence isomorphic to
$C_3$ if $q-1$ is a multiple of $3$ and trivial otherwise.
Let 
$$d := |Z| = \gcd(3,q-1) \mbox{ and } 
\omega := \exp \left(\frac{2\pi i}{3} \right) .$$ 

\subsection{The characters of $\GL_2(q)$}

As we use Harish-Chandra induction from the Levi factor 
$\GL_2(q)$ of the subgroup $P$, we first deal with 
this group.

\begin{theorem} \label{GL2}
	Assume that $q$ is a power of some odd prime $p$. Then
	all characters $\psi \in \Irr ^{o}(\GL_2(q))$ have unitary discriminant 
	$\disc(\psi ) = (-1)^{\psi(1)/2} $.
\end{theorem}

\begin{proof}
	If $\psi $ is not faithful, then the character field of $\psi $ 
	is real (see \cite{Steinberg}). 
	So $\psi $ is a faithful irreducible character of $\GL_2(q)$.
	Then the restriction of $\psi $ to 
	the subgroup $D_8 \leq \GL_2(q) $ is a multiple of 
	the unique character of $D_8$ that restricts non-trivially 
	to the center of $D_8$. This character has degree 2, indicator +, 
	trivial Schur indices and orthogonal discriminant $-1$.
	As in Theorem \ref{extraspecial} we conclude that 
	the unitary discriminant is $\disc(\psi ) = (-1)^{\psi(1)/2} $.
\end{proof}

\begin{remark} \label{GL2even}
        Assume that $q$ is a power of $2$. Then 
	$\GL_2(q)$ has $q-1$ irreducible characters of degree $q$. 
	All these characters restrict to the Steinberg character 
	of degree $q$ of $\SL_2(q)$ and hence have 
	unitary/orthogonal discriminant $(-1)^{q/2} (q+1)$ (see \cite[Theorem 6.2]{SL2}).
\end{remark}

\subsection{The characters of $P$}

Now let $\chi  \in \Irr(P)$ be an irreducible character of $P\cong \mathbb{F}_q^2 \rtimes \GL_2(q) $ restricting 
non-trivially to the 
 abelian normal subgroup $A:= \mathbb{F}_q^2 $ of order $q^2$. 
 As $P$ acts transitively on $A\setminus \{ 1 \}$,
 the restriction of $\chi $ to $A$ is a multiple of 
 the sum of all non-trivial linear characters $\psi $ of $A$. 
 By Clifford theory, the character $\chi $ is induced up from a character of 
 the inertia subgroup $T_{\psi } \cong  H\ltimes A$ of any of these characters $\psi $, 
 where 
 $$H = \left\{ \left( \begin{array}{cc} a & b \\ 0 & a^{-2} \end{array} \right) 
	 \mid  a \in \mathbb{F}_q^{\times }, b\in \mathbb{F}_q \right\} \leq \GL_2(q) $$
	 is isomorphic to $(\mathbb{F}_q^{\times } \ltimes  \mathbb{F}_q)   $. 
	 The 
 center of $T_{\psi }$ is $Z$  and the 
 index of $T_{\psi }$ in $P$ is $q^2-1$. 
The irreducible characters $\phi $ of $H$ consist of $(q-1) $ linear characters 
and $d^2$ characters of degree $(q-1)/d$. 
Inducing the characters $\phi \otimes \psi $ with $\phi \in \Irr(H)$ 
from $H\ltimes A$ to $P$ we obtain $(q-1)$ irreducible characters 
of degree $(q^2-1)$ and $d^2$ characters of degree $(q-1)^2(q+1)/d$.  

\begin{remark}\label{osP}
	The irreducible characters $\chi \in \Irr(P)$ for which 
	the restriction of $\chi $ to a Sylow $p$-subgroup of $P$ does not 
	contain the trivial character are exactly the $d^2$ characters 
	of degree $\chi (1) = (q-1)^2(q+1)/d$. 
	By Theorem \ref{Benard} all these characters have trivial Schur index. 
	The $d$ characters $\chi $ that restrict trivially to the center 
	are rational, the other $d^2-d$ characters $\chi $ have character field
	$\mathbb{Q}(\omega )$.
	All these characters have trivial unitary, resp. orthogonal, discriminant.
\end{remark}

\begin{proof}
Write $q-1= a b$ such that $a$ is a power of $3$ and
  $b$ is not a multiple of $3$ and put
        $$U:=\langle P' , g^a \mid g\in P \rangle $$ to be the normal subgroup of
        index $a$ in $P$ (so $U=P$ if $q\not\equiv 1 \pmod{3} $).
	Let $\chi \in \Irr(P)$ be one of the irreducible characters of 
	degree $(q-1)^2(q+1)/d$ from Remark \ref{osP}. 
	By Clifford theory the restriction $\chi _{|U }$ of $\chi $ to $U$ is 
	a sum of a 3-power number of conjugate characters of the same degree, 
	in particular the degrees of the constituents of $\chi _{|U}$ are even. 
	Moreover these constituents are rational and hence  $\chi _{|U }$ 
	is an orthogonally stable rational character
        that restricts orthogonally stably to a Sylow $p$-subgroup
	of $U$. Now \cite[Theorem 4.3 and Corollary 4.4]{Albanian} yields 
	$$\disc (\chi _{|U}) = (-1)^{\chi(1)/2} p^{\chi(1)/(p-1)}  (\mathbb{Q}^{\times })^2 = 1 $$
        and hence also the unitary discriminant of $\chi $ is $1$.
\end{proof}

\subsection{The characters of $\SL_3(q)$}

\begin{theorem} \label{SL3A}
	If $q$ is odd then all characters 
	$\chi \in \Irr ^{o} (\SL_3(q) )$ have unitary discriminant $(-1)^{\chi(1)/2}$.
\end{theorem}

\begin{proof}
	We use the description in \cite{Steinberg} of the characters of 
	$\SL_3(q)$. 
	The characters $\chi $ of degree $(q+1)(q^2+q+1)$ and 
	$(q-1)(q^2+q+1)$ are Harish-Chandra induced from 
	the maximal parabolic subgroup 
	$P$ from a 
	a character $\psi $ of 
	degree $(q+1)$ resp. $(q-1)$ of the Levi complement $\GL_2(q)$. 

	If $\chi(1) = (q-1)(q^2+q+1)$,
	then there is a unique such character $\psi $ with $\chi = \psi^G$. 
	In particular  $\mathbb{Q}(\chi ) = \mathbb{Q}(\psi )$.
If $\mathbb{Q}(\psi ^G)$ is not real, then also $\psi \in \Irr^{o}(\GL_2(q))$
	and Theorem \ref{GL2} shows that the 
	unitary discriminant of $\psi $ 
	and hence the one of $\chi = \psi^G$ is $(-1)^{(q-1)/2}$.

	If $\chi(1) = (q+1) (q^2+q+1)$, then there are three such characters 
	$\psi $ inducing to the same character $\chi $. 
	As the Galois group of 
	 $\mathbb{Q}(\psi ) / \mathbb{Q} (\chi )$ 
	 acts on the constituents of degree $q+1$ of $\chi _{|P}$
	 it has at least one orbit of odd length. 
	 So one of these characters $\psi $ 
	 satisfies $[\mathbb{Q}(\psi ) : \mathbb{Q} (\chi )]$ is odd.
	 By Lemma \ref{odd}, this
	again allows to conclude that the unitary discriminant
	of $\chi $ is the same as the one of $\psi $.

	It remains to consider the characters $\chi $ of degree $(q-1)^2(q+1)$ 
	and $(q-1)^2(q+1)/3$ (if $q\equiv 1 \pmod{3}$).  

	For both degrees, the character $\chi $ does not appear 
	in the permutation character of $G$ on the 
	$(q-1)^2(q+1)$ cosets of a Sylow $p$-subgroup $S$ of $G$: 
	This is clear if $\chi(1) = (q-1)^2(q+1) = [G:S]$ because 
	all constituents of $1_S^G$ have degree $\leq [G:S]-1$. 
	If $\chi (1) = (q-1)^2(q+1)/3$, we note that the 
	center  of $G$ has order 3 and orbits of length 
	3 on the cosets of $G/S$.
	So if $\chi $ occurs in $1_S^G$, then this permutation character 
	has three distinct constituents of degree $(q-1)^2(q+1)/3$ 
	leading to the same contradiction as before. 

	In particular, the restriction of $\chi $ to $P$ is a sum of the
	characters from Remark \ref{osP}, showing again that the unitary 
	discriminant of $\chi $ is trivial.

	From the tables in \cite{Simpson1973TheCT},
	we see that the characters of degree $q(q+1)$ and 
	$(q+1)(q^2+q+1)/3$ are rational. Their orthogonal discriminant
	can be read off from \cite[Theorem 4.7]{SL3SU3}.  
\end{proof}

\begin{theorem} \label{SL3B}
	 If $q$ is even then all characters 
	$\chi \in \Irr ^{o} (\SL_3(q) )$ have degree $q(q^2+q+1)$ and 
	unitary discriminant $(-1)^q (q+1)$.
\end{theorem} 

\begin{proof}
	For 2-powers $q$, the even degree 
	irreducible characters of $\SL_3(q)$ are as follows:
	\begin{itemize}
		\item[(i)] one character of degree $q^2+q$ and Frobenius-Schur indicator 
	$+$, 
\item[(ii)] the Steinberg character of degree $q^3$  and Frobenius-Schur indicator $+$, 
\item[(iii)] 
	 $q-1$ characters of degree $q(q^2+q+1)$, one of which is rational and
	 of indicator $+$, while the others have indicator 'o' (see \cite{Steinberg}, 
	 \cite{Simpson1973TheCT}). 
	\end{itemize}
	From \cite[Table VI]{Steinberg} and 
	 the arguments given just before we conclude that the 
	 $q-1$ characters from (iii) 
	 are Harish-Chandra induced from the $q-1$ irreducible characters $\psi $ of 
	 degree $q$ of the Levi factor $\GL_2(q)$ of $P$. 
	 It follows that 
	 they have the same character field ${\mathbb Q}(\chi ) = {\mathbb Q}(\psi )$ and discriminant. From Remark \ref{GL2even} we get $\disc(\psi ) = (-1)^{q/2}(q+1)$ whenever 
	 the character field $\Q(\chi )$ is not real, i.e. $\chi \in \Irr^o(\SL_3(q))$.
	 As the index of $P$ in $\SL_3(q)$ is odd,
	 also the discriminant of $\chi $ is $(-1)^{q/2} (q+1)$.
\end{proof}
 
\section{The special unitary group, general results}  \label{SU3gen}

\subsection{The special unitary group}  \label{SU3}

Let $p$ be a prime and let $q$ be a power of $p$.
The special unitary group $\mathrm{SU}_3(q)$ is the
stabiliser in $\mathrm{SL}_3(q^2)$ of a non-degenerate Hermitian form
on $\mathbb{F}_{q^2}^3$.
Up to isometry there is a unique such form.
We put
$$\Omega= \begin{pmatrix}
0 & 0 & 1 \\
0 & 1 & 0 \\
1 & 0 & 0
\end{pmatrix}$$ and denote by $\Phi $ the $\mathbb{F}_q$-linear map on
$\mathbb{F}_{q^2}^{3\times 3} $ that raises each matrix entry to
the $q$-th power.
Then
$$
    \mathrm{SU}_3(q) :=  \{ \g \in \mathrm{SL}_3(q^2) \mid \Phi(\g)^{tr} \cdot \Omega \cdot \g = \Omega \}.
$$
Let
$$
B \coloneqq \left\{
\begin{pmatrix}
d & a & b \\
0 & e & c \\
0 & 0 & f
\end{pmatrix} \in \SU_3(q)
\right\} \mbox{ and }
U \coloneqq \left\{
\begin{pmatrix}
1 & a & b \\
0 & 1 & c \\
0 & 0 & 1
\end{pmatrix} \in \SU_3(q)
\right\}.
$$
Then $U$ is the unipotent radical of $B$ and a Sylow $p$-subgroup of $\SU_3(q)$, and $B=N_{\SU_3(q)}(U) = U \rtimes T$
is a (standard) Borel subgroup,
where $T:=\{\mathrm{diag}(d,e,f)\in \SU_3(q)\}$ is a maximal torus. 

We also put $$\w := \begin{pmatrix}
0 & 0 & 1 \\
0 & -1 & 0 \\
1 & 0 & 0
\end{pmatrix} $$ 
to denote a generator of the Weyl group of $\SU_3(q)$. 
Then $$\SU_3(q) =  B \disj B\w B \mbox{ and } B\cap \w B \w = T .$$

We fix an element $\T \in T$ such that $T = \langle \T \rangle \cong C_{q^2-1} $ 
and put 
$$\T_0 := \T^{(q^2-1)/2} =  \diag (-1,1,-1) $$ 
to denote the element of order 2 in $T$.

The center $Z  := Z(U) \cong (\mathbb{F}_q,+)$ of the unipotent radical is generated as a normal 
subgroup of $B$ by 
$$\z:=
\begin{pmatrix}
1 & 0 & b \\
0 & 1 & 0 \\
0 & 0 & 1
\end{pmatrix} \in Z,$$ where $b \in \mathbb{F}_{q^2}^{\times}$ is an element such that $b+b^q=0.$
Then $\T_0$ commutes with the elements of $Z$ and inverts the classes in $U/Z$, i.e. 
$$\T_0 \h Z \T_0  = \h^{-1} Z $$ for all $\h \in U$.
We denote by $A_0$ the subgroup 
$$A_0 := \langle \z , \T \rangle \cong (\F_q,+) \rtimes (\F_{q^2}^{\times },\cdot ) .$$
We also need two rational quaternion algebras: 
$${\mathcal Q}_p := {\mathcal Q}_{\Q}(\infty, p) \mbox{ and } 
{\mathcal Q}_2 := {\mathcal Q}_{\Q}(\infty, 2) = (-1,-1)_{\Q }  $$ 
ramified at the infinite place and the prime $p$ resp. $2$.
For $p=2$ we have that ${\mathcal Q}_2 = {\mathcal Q}_p$ and 
for $p\equiv 3 \pmod{4} $ the algebra ${\mathcal Q}_p = (-1,-p)_{\Q } $.

\subsection{Restrictions to the
Borel subgroup $B$}  

The character table and the conjugacy classes of $B$ have already been obtained in 
\cite[Table 2.1]{Geck}. That paper also contains the character table of $\SU_3(q) $
using the same notation as in  \cite{Simpson1973TheCT}. 
The index of the character names always gives their degree. 
To give the characters of $B$ in a unified way we denote by
	$$d:=\gcd (3,q+1) \in \{1,3\}$$ 
	 the order of the center of $\SU_3(q)$, which is also the center of $B$. 
	 We put 
	 $$\delta := \exp(2\pi i/(q^2-1))  \mbox{ and } \omega := \exp(2\pi i/3) $$
	 to denote a primitive $(q^2-1)$th, resp. third, complex root of unity.

	 \begin{proposition} \cite[Table 2.1 b)]{Geck}  \label{IrrB}
The irreducible characters of $B$ are as follows: 
 \begin{itemize} 
	 \item[a)] The  $q^2-1$ linear characters  $\vartheta _1^{(u)} $ 
		 for $0\leq u \leq q^2-2$ 
		 with character field $\Q( \delta ^u )$.
	 \item[b)] The $q+1$ characters $\vartheta _{q^2-q}^{(u)} $ 
		 for $0\leq u \leq q$ with character field 
		 $\Q(\delta ^{(q-1) u})$.
		 Then $[\vartheta _{q^2-q}^{( 0 )}] = [{\mathcal Q}_p ]$ and
		 all the other characters have trivial Brauer element.
	 \item[c)] The $d^2$ characters $\vartheta _{(q^2-1)/d}^{(u,v)}$ 
		 for $0\leq u, v \leq d-1$ with character field 
		 $\Q(\omega^u )$.
 \end{itemize}
	 \end{proposition}

		 Note that for $d=1$ the character 
		 $\vartheta _{(q^2-1)/d}^{(0,0)} $ is denoted by
		 $\vartheta _{q^2-1} $ in \cite{Geck}.

We want to describe the restrictions of the irreducible characters of $\SU_3(q)$ to $B$. For that, we introduce two further families of characters of $B$.

\begin{definition} \label{phiundtau}
We define
\[ \varphi^{(v)}:= \left( d\sum_{j=1}^{(q+1)/d} \vartheta_{q^2-q}^{(dj+v)} \right) , \ \tau^{(u)}=\sum_{j=0}^{d-1}  \vartheta_{(q^2-1)/d}^{(u,j)}     \]
for arbitrary integers $v$ and $u$. 
Note that if $d=1$, then $\tau^{(0)}=\vartheta_{q^2-1}^{(0,0)}$. 
Here and in the following we consider the upper index $u$ of 
	$\vartheta _1^{(u)}$ modulo $q^2-1$ and the one of 
	$\vartheta _{q^2-q}^{(u)}$ modulo $q+1$. 
\end{definition}

\begin{proposition}\label{restoB}
	The restrictions to $B$ of the characters of $\SU_3(q)$ 
 are given by:
\begin{enumerate}[label=\roman*)]
	\item $\mathrm{Res}_B(\chi_1)=\vartheta_1^{(0)}$.
	\item $	\mathrm{Res}_B(\chi_{q(q-1)})=\vartheta_{q^2-q}^{(0)}$.
\item $\mathrm{Res}_B(\chi_{q^3})=\vartheta_1^{(0)}+\tau^{(0)}+\varphi^{(0)}-\vartheta_{q^2-q}^{(0)}$.
    \item $\mathrm{Res}_B(\chi_{q^2-q+1}^{(u)})=\vartheta_1^{((q-1)u)}+\vartheta_{q^2-q}^{(u)}$, for $1\leq u \leq q$.
    \item   $\mathrm{Res}_B(\chi_{q(q^2-q+1)}^{(u)})=\vartheta_1^{((q-1)u)}+ \tau^{(u)} + \varphi^{(u)} - \vartheta_{q^2-q}^{(u)} -\vartheta_{q^2-q}^{(-2u)}$, for $1\leq u\leq q$. 
    \item   $\mathrm{Res}_B(\chi_{(q-1)(q^2-q+1)}^{(u,v,w)})= \tau^{(u+v+w)} + \varphi^{(u+v+w)} - \vartheta_{q^2-q}^{(u-2v-2w)}- \vartheta_{q^2-q}^{(v-2u-2w)}- \vartheta_{q^2-q}^{(w-2u-2v)}$, \\ for
	    $1\leq u < v \leq (q+1)/d, v<w\leq (q+1), u+v+w \equiv 0 \pmod{q+1} $. 
     \item $\mathrm{Res}_B(\chi_{(q-1)(q^2-q+1)/3}^{(u)})= \vartheta_{(q^2-1)/3}^{(0,u)} + 1/3\varphi^{(0)} - \vartheta_{q^2-q}^{(0)}$, for $0\leq u \leq 2$. 
	\item $\mathrm{Res}_B(\chi_{q^3+1}^{(u)})= \vartheta_1^{(u)}+\vartheta_1^{(-qu)} +\tau^{(u)} + \varphi^{(u)} - \vartheta_{q^2-q}^{(u)}$, \\
		for $1\leq u \leq q^2-1, (q-1) \nmid u.$ Note that $\chi_{q^3+1}^{(u)}=\chi_{q^3+1}^{(-uq)}$.
\item $\mathrm{Res}_B(\chi_{(q+1)(q^2-1)}^{(u)})= \tau^{(u)} + \varphi^{(u)}$, for $1\leq u < q^2-q+1,  (q^2-q+1)/d \nmid u .$ \\
Note that $\chi_{(q+1)(q^2-1)}^{(u)}=\chi_{(q+1)(q^2-1)}^{(-uq)}=\chi_{(q+1)(q^2-1)}^{(uq^2)}$ where the upper indices are understood modulo $q^2-q+1$.
    \item $\mathrm{Res}_B(\chi_{(q+1)(q^2-1)/3}^{(u)})= \vartheta^{(1,u)}_{(q^2-1)/3} + 1/3\varphi^{(1)}$, for $0\leq u \leq 2$. 
    \item $\mathrm{Res}_B(\chi_{(q+1)(q^2-1)/3}^{(u)'})= \vartheta^{(2,u)}_{(q^2-1)/3} + 1/3\varphi^{(2)}$, for $0\leq u \leq 2$. 
\end{enumerate}
Note that the characters $\chi_{(q-1)(q^2-q+1)/3}^{(u)}$, $\chi_{(q+1)(q^2-1)/3}^{(u)}$ and $\chi_{(q+1)(q^2-1)/3}^{(u)'}$ only exist if $d=3$.
\end{proposition}

\section{The unitary discriminants of $\SU_3(q)$ for $q$ even} \label{SU3even}

In this section we assume that  $q\geq 4$ is a power of $2$.
As $\chi_{q(q-1)}$ and $\chi_{q^3}$ are rational characters we get that
$$\Irr^o(\SU_3(q)) = \{ \chi_{q(q^2-q+1)}^{(u)} \mid 1\leq u \leq q \} .$$ 

The aim of this section is the proof of the following theorem. 

\begin{theorem} \label{mainSU3even}
	Let $q\geq 4$ be a power of $2$. 
	Then $\disc(\chi ) = q+1$ for all 
 $\chi \in \Irr^{o}(\SU_3(q))$.
\end{theorem}

From Proposition \ref{restoB} we get 
\begin{remark} \label{defdet}
	For $1\leq u\leq q$ the sum
	$$M^{(u)} := \chi_{q(q^2-q+1)}^{(u)} + \chi_{(q^2-q+1)}^{(u)} = (\vartheta _1^{((q-1)u)}) ^{\SU_3(q)} $$
	is the monomial character induced from the character $\vartheta _1^{((q-1)u)}$
	of $B$. 
	The character fields of $M^{(u)}$, $\chi_{q{(q^2-q+1)}}^{(u)}$, $\chi_{(q^2-q+1)}^{(u)}$ and 
	$\vartheta_1^{((q-1)u)}$ are identical and isomorphic to $K^{(u)}:= {\mathbb Q}(\delta ^{(q-1)u})$.
	Let $V^{(u)}$ denote the $K^{(u)} \SU_3(q) $ module affording the 
	character $M^{(u)}$.
	Then $V^{(u)}$ is an orthogonal sum of the 
	two absolutely irreducible submodules $V_{q{(q^2-q+1)}}^{(u)}$ and $V_{(q^2-q+1)}^{(u)}$. 
	As  $\SU_3(q)$ fixes the standard form $I_{q^3+1}$ on $V^{(u)}$
	the discriminant of $\chi_{q{(q^2-q+1)}}^{(u)} $ 
	is 
the discriminant of the restriction of the 
standard form to the submodule 
$V_{(q^2-q+1)}^{(u)}$.
\end{remark} 

\begin{remark} \label{basismon}
	In the notation of Section \ref{SU3}, we have
		$$\SU_3(q) = B \cup B\w B = B \disj \Disj _{\h \in U} B \w \h .$$
	So a basis of  $V^{(u)}$
	is given by $\{ B \} \cup \{ B\w \h \mid \h \in U \} $. 
\end{remark}

In this notation we obtain the following 
\begin{lemma}\label{schur}
	The Schur basis of $\End _{\SU_3(q)} (V^{(u)})$ is 
	$(I_{q^3+1} , E )$ where 
	$$E_{B,B} = 0, E_{B,B\w \h} = 1 , E_{B\w \h,B} = 1 \mbox{ for all } \h \in U.$$ 
\end{lemma}

\begin{proof}
	By the well known formulas for the Schur basis elements 
	(see for instance \cite[Proposition (1.10)]{HabilJuergen}) 
	we have  $E_{B,B}=0$  and 
	$E_{B,B\w \h} = \vartheta_1^{((q-1)u)} (\h) = 1  $. 
To compute 
	$E_{B\w \h,B} $ we need to write $\h^{-1} \w = \g \w \h' $ for $\h'\in U$, $\g \in B$. 
	But then $\g^{-1} \h^{-1} = \w \h' \w^{-1} $ is an element of $2$-power order
	in $B$ and hence also $\g \in U$, so $\vartheta_1^{((q-1)u)}  (\g) = 1$, 
	whence 
	$$E_{B\w \h,B} = \vartheta_1^{((q-1)u)} (\g) \vartheta_1^{((q-1)u)} (\h') = 1 .$$
\end{proof}

\begin{lemma}\label{eigenvals} 
	The eigenvalues of $E$ are 
	$ \epsilon  q$ and $ -\epsilon q^2  $ for some $\epsilon \in \{ 1,-1 \} $ with multiplicities 
	$q(q^2-q+1)$ and $(q^2-q+1) $.
\end{lemma} 

\begin{proof}
	The 2-dimensional $K^{(u)}$-space generated by $(I_{q^3+1}, E )$ is 
	a ring, in particular $E^2$ is a $K^{(u)}$-linear combination of 
	$E$ and $I_{q^3+1}$ and $E$ has exactly 2 distinct eigenvalues.
	Also the multiplicity of the eigenvalues are given by the 
	dimensions of the irreducible constituents of $M^{(u)}$.
	Let $a$ denote the eigenvalue of $E$ occurring with multiplicity 
	$q^2-q+1$ and $b$ the one with multiplicity $q(q^2-q+1)$. 
	Then 
	$$ a(q^2-q+1) + bq(q^2-q+1) = \mbox{trace}(E) = 0  ,$$
	 so  
	$a=-bq .$
	As the diagonal entries of $E$ are 0 and the first diagonal entry of 
	$E^2$ is $q^3$ by Lemma \ref{schur}, we obtain 
	the constant coefficient of the minimal polynomial of $E$ as
	$$-q^3 = ab= -b^2q \mbox{ and hence } b = \epsilon  q $$
	for some $\epsilon  = \pm 1 $.
\end{proof}

\begin{corollary} 
	The rows of $E+\epsilon  q^2I_{q^3+1} $ span the submodule 
	$V_{q(q^2-q+1)}^{(u)} $ and 
	the rows of $E-\epsilon  q I_{q^3+1}$ the submodule 
	$V_{(q^2-q+1)}^{(u)}$.
\end{corollary} 

As $U$ is a normal subgroup of $B$, we obtain a basis of the 
$2$-dimensional $B$-eigenspace as follows. 

 \begin{remark} \label{Wu}
	 The $B$-eigenspace for the character $\vartheta _1^{((q-1)u)} $ in $V^{(u)}$ is 
 $ W^{(u)}:=\langle B , \sum _{\h \in U} B\w \h   \rangle  $.
	Put 
	 $$v^{(u)}:=-\epsilon  q^2  B + \sum _{\h \in U} B\w \h   \mbox{ and }
	 w^{(u)}:=\epsilon  q B + \sum _{\h \in U} B\w \h  .$$
	 where $\epsilon  $ is as in Lemma \ref{eigenvals}. 
	 Then
	 $$\langle v^{(u)} \rangle = W^{(u)} \cap V_{q(q^2-q+1)}^{(u)} \mbox{ and } 
	 \langle w^{(u)} \rangle = W^{(u)} \cap V_{(q^2-q+1)}^{(u)} .$$
 \end{remark}

 The next lemma completes the proof of Theorem \ref{mainSU3even}.

 \begin{lemma} 
	 The unitary discriminant of $\chi_{q(q^2-q+1)}^{(u)}$ is represented by $q+1$.
 \end{lemma} 

 \begin{proof}
	The restriction of $M^{(u)}$ to $B$ is 
	 $$ \mathrm{Res}_B(M^{(u)}) =  \vartheta_1^{((q-1)u)}+\vartheta_{q^2-q}^{(u)} +  \mathrm{Res}_B(\chi_{q(q^2-q+1)}^{(u)}) .$$
	 The restriction of the character
                 $\vartheta_{q^2-q}^{(u)}$ to the Sylow 2-subgroup $U$ of $B$
                        is a sum of $q-1$ rational characters of degree $q$.
	 As $q$ is a multiple of 4, Corollary 
	 \ref{2groups} implies that this restriction  
 is unitary stable of unitary discriminant $1$.
So the discriminant of the submodule
	 $V_{(q^2-q+1)}^{(u)} $ of $V^{(u)}$ is 
	 the square length of $w^{(u)}$,  which is 
	 $q^2+q^3=(q+1)q^2$. 
	 As the product of the discriminants of the two submodules 
	 $V_{(q^2-q+1)}^{(u)} $ and $V_{q(q^2-q+1)}^{(u)} $ of $V^{(u)}$ is $\disc(V^{(u)}) = 1$, 
	 we obtain 
	 $$\disc (V_{q(q^2-q+1)}^{(u)} ) = q+1  = \disc(\chi_{q(q^2-q+1)}^{(u)} ) .$$
 \end{proof}

\section{The unitary discriminants of $\SU_3(q)$ for $q$ odd} \label{SU3odd}

\subsection{The strategy for $q$ odd} 
Let $q$ be odd and let $\chi \in \Irr ^o (\SU_3(q))$ be an even degree indicator 'o' irreducible 
complex character. The goal of this subsection is to outline the strategy which we will use in the sequel to calculate the unitary discriminant of $\chi $.

The restriction of $\chi $ to $B$ decomposes as 
$$\mathrm{Res}_B(\chi ) = \chi_T + \chi _U ,$$ where 
$\chi _T$ is the $U$-fixed part of $\mathrm{Res}_B(\chi)$. 
The character $\chi _T$ is also known as the 
\textit{Harish-Chandra restriction} of $\chi$.

For orthogonal characters $\chi $ the restriction
$R:=\mathrm{Res}_U(\chi _U)$ to $U$ is an orthogonally stable character 
of the $p$-group $U$, so 
the determinant of $R$ and hence the orthogonal determinant of $\chi_U$ 
is given in \cite[Theorem 4.3 and Corollary 4.4]{Albanian}. 
As $q$ is odd, $R$ is a sum of 
odd degree characters and hence never unitary stable. 
However, it turns out that $\chi _U$ is a unitary stable character 
of $B$ and hence has a well defined unitary discriminant. 
To compute $\disc(\chi_U)$, we restrict further to 
the metabelian subgroup $A_0\leq B$ from Section \ref{SU3} (Section \ref{Aunitary}).

\begin{remark} \label{monirr} 
The complex irreducible characters $\chi $ of $\SU_3(q)$ for which $\chi _T$ is non-zero 
are monomial characters on the set of $\chi(1) = q^3+1$ isotropic points in the 
natural 3-dimensional unitary geometry. 
Here we use condensation techniques to find $\disc(\chi_T)$: 
\\
Put $J_U:= \frac{1}{|U|} \sum _{\h \in U} \h $ to denote the projection 
onto the $U$-fixed space and let $V_{\chi }$ be a $\Q(\chi ) \SU_3(q)$-module 
affording the character $\chi $. 
Then the decomposition $\chi  = \chi_T + \chi_U$ corresponds to the orthogonal 
decomposition
$$V_{\chi} = V_{\chi } J_U \perp V_{\chi }(1-J_U ) $$ 
and hence $$\disc(\chi) = \disc(V_{\chi} J_U) \disc(\chi_U) .$$ 
The discriminant of the 2-dimensional module $V_{\chi} J_U$
is computed by obtaining the action of 
$J_U \T J_U $ and $J_U \w J_U $ on this module. 
\end{remark}

\subsection{The unitary discriminants of $A_0$}  \label{Aunitary}

In this section we compute the unitary discriminants of 
the subgroup $A_0$ of $B$ defined in Section \ref{SU3},
which are then used in Section \ref{Bunitary} to
obtain the unitary discriminants of the irreducible characters of 
degree $q^2-q$ of $B$.

Recall that 
$A_0 = \langle \z , \T \rangle \cong (\F_q,+) \rtimes (\F_{q^2}^{\times },\cdot ) .$
The subgroup
$\langle Z(U) , \T^{q-1} \rangle $ is an abelian normal subgroup of 
order $q(q+1)$ and index $(q-1)$ in $A_0$. 
The center  $Z(A_0) \cong C_{q+1}$ is generated by $\T^{q-1}$. 
Let 
$\lambda ^{(u)} $ be the linear character of $Z(A_0)$ defined by 
$$\lambda ^{(u)} (\T^{q-1}) = \delta ^{(q-1) u} \mbox{ for } u=0,\ldots, q.$$ 
The group $A_0$ has $(q^2-1)$ linear characters, the ones that 
restrict trivially to $Z(U)$, and 
$q+1$ irreducible characters of degree $q-1$, $\mu _{q-1}^{(u)} $ for $u=0,\ldots, q$,
that restrict to $Z(A_0)$ as $(q-1) \lambda ^{(u)} $.

\begin{theorem}\label{A0deg(q-1)}
	The non-linear irreducible characters of $A_0$ are the
	characters $\mu_{q-1}^{(u)}$, $0\leq u\leq q$ 
	of degree $q-1$ and character field 
	$\Q(\mu_{q-1}^{(u)}) = \Q(\delta^{(q-1) u}) .$
        \begin{itemize}
            \item[a)] $\mu _{q-1}^{(0)}$ is the character of a rational representation. Its orthogonal discriminant is $(-1)^{(q-1)/2}q$.
	\item[b)] $\mu _{q-1}^{((q+1)/2)}$
is a rational character of Frobenius-Schur indicator $-$. Its
		Brauer element is $[\mu_{q-1}^{((q+1)/2)}] = [{\mathcal Q}_p]$.
                \item[c)]
			For $u \not\in \{ 0, (q+1)/2 \} $ the Frobenius-Schur 
			indicator of $\mu_{q-1}^{(u)}$ is \textrm{'o'} and 
                        the discriminant is
			$$\disc(\mu_{q-1}^{(u)}) = (-1)^{(q-1)/2} q^{u-1} = 
			\left\{ \begin{array}{ll}
		-q & \mbox{ if }  (u,q) \equiv (0,3) \pmod{(2,4)} \\
		-1 & \mbox{ if } (u,q) \equiv (1,3) \pmod{(2,4)} \\
		q & \mbox{ if }  (u,q) \equiv (0,1) \pmod{(2,4)} \\
		1 & \mbox{ if } (u,q) \equiv (1,1) \pmod{(2,4)}. \\
                        \end{array} \right. $$
        \end{itemize}
\end{theorem}

To prove the theorem, we restrict further to a subgroup $A$ of $A_0$:
Write $q+1 =  e b$, where $b$ is odd and $e$ is a power of $2$ and 
put $\T_1 := \T^b $ a generator of the subgroup of order $e(q-1)$ of 
the torus $T$. 
Put 
$$A := \langle \z , \T_1 \rangle  .$$
Then $A_0 = A \times \langle \T^{(q-1)e} \rangle $ and the 
irreducible characters of $A_0$ are obtained as tensor product 
$\mu \otimes \lambda $, where $\mu $ is an irreducible character of $A$ and
$\lambda $ a linear character of the cyclic group $ \langle \T^{(q-1)e} \rangle $ 
of order $b$.

We first compute the unitary/orthogonal discriminants and Schur indices of 
the even degree irreducible characters of $A$.

\begin{lemma}\label{Adeg(q-1)}
	For $u=0,\ldots , q$ we put
	$\mu := \Res_A(\mu_{q-1}^{(u)})$.
	Then $\mu $ only depends on $u\pmod{e}$.
	\begin{itemize}
\item[a)] If $e$ divides $u$, then $\mu $ is the character of a rational representation. 
Its orthogonal discriminant is $(-1)^{(q-1)/2}q$. 
\item[b)] If $u \pmod{e} = e/2$, then $\mu $ 
is a rational character of Frobenius-Schur indicator $-$. Its 
			Brauer element is $[\mu] = [{\mathcal Q}_p ]$,
			the class of the rational quaternion algebra 
			ramified only at  $p$ and $\infty $.
		\item[c)]
			If $u \pmod{e} \not\in \{ 0, e/2 \} $, then $e \geq 4$ and
			hence $q\equiv 3 \pmod{4}$ is an odd power of a prime
			$p\equiv 3 \pmod{4}$. 
			Then the character field  of 
			$\mu $ is 
			$\Q(\mu) = \Q (\zeta_e^u) $ and contains a primitive
			fourth root of unity.
			The unitary discriminant of $\mu $ is
	$$\disc(\mu) = \left\{ \begin{array}{ll} 
				-q & \mbox{ if } u \mbox{ is even} \\
				-1 & \mbox{ if } u \mbox{ is odd.} 
			\end{array} \right.
			$$
	\end{itemize}
\end{lemma}

\begin{proof}
	Put $Z:=Z(U)$ and $\T_2 := \T^{(q-1)b} = \T_1^{q-1}$  to denote a generator of 
	the center of $A$.
Then the character $\mu  $ is a 
monomial character induced from a 
linear character, $\lambda $,
	of the normal subgroup $Z \times \langle \T_2 \rangle $. 
	Theorem \ref{Yamada}
	says that, for any field $K$, the Brauer element $[\mu ]_{K(\mu)}$ of
$\mu$ is the inverse of the class of 
	the crossed product algebra ${\mathcal Q}:=(K(\lambda ) , \Gamma )$ in 
	the Brauer group of $K(\mu )$,
where $\Gamma $ is the Galois group of $K(\lambda )/K(\mu )$. 
	In our case, $K(\lambda )$ is generated by 
	$\lambda (Z) =\langle \exp(2\pi i/p ) \rangle $ and 
	$\lambda (\T_2)   = \zeta_e^u $.
	Let $e_0 := e/\gcd(e,u) $ denote the order of $\lambda (\T_2) $.
	For $K={\mathbb{Q}} $ we have $\Gamma = C_{p-1} = \langle \sigma \rangle = \Gal ({\mathbb{Q}}(\exp(2\pi i /p))/{\mathbb{Q}}) $, where the cocycle 
	is given by $\sigma ^{p-1} =\zeta_e^u  $. 

	Let 
	$\langle a \rangle = ({\mathbb{Z}}/p{\mathbb{Z}})^{\times }$ and put
$$A':= \langle \z' , \T_2' , \sigma  \mid 
	\z'^p = 1, \T _2'^{e_0} = 1, \z'^{\sigma } = \z'^a , \sigma ^{p-1} = \T_2' \rangle. $$ 
	Then $|A'| = p (p-1)e_0$ and $A'$ is a group all of whose 
	Sylow subgroups are cyclic. 
	Moreover, $[{\mathcal Q}] = [\chi ]$ for an 
	irreducible faithful character $\chi $ of degree $p-1$
	of the group $A'$.  

	\begin{itemize}
		\item[a)] If $e_0 = 1$  and $K={\mathbb{Q}} $, then 
			${\mathcal Q}={\mathbb{Q}}^{p-1\times p-1}$, so we get the 
			Schur indices in a). 
			For the orthogonal determinant, note that 
			the representation corresponding to 
			$\mu$ fixes the root lattice $A_{q-1}$ of determinant $q$.
		\item[b)] If $e_0 = 2$, so $\lambda (\T_2) = -1$, then 
			again $\chi $ is rational. Now the Frobenius-Schur indicator
			of $\chi $ is $-$. To compute the local Schur indices of 
			$\chi $ we use Theorem \ref{Benard}. 
			Over the completions at primes dividing $p-1$ this Schur-index is 1, so it remains to compute $m_{\Q_p}(\chi )$. 
			Any $p$-modular constituent of $\chi $ is a (faithful) representation of $A'/\langle \z' \rangle $, a cyclic group of order $(p-1)e_0 = 2(p-1)$ and hence
			it character field is the cyclotomic field of order $2(p-1)$ 
			over $\Q_p$. This has degree 2 over $\Q_p$, and hence 
			$m_{\Q_p}(\chi) =2$. 
			So $[\chi ] = [{\mathcal Q}_p]$ is the class of the
			rational quaternion algebra ramified at 
			$p$ and $\infty $.
		\item[c)] If $e_0 \geq 4$, then the indicator of $\chi $ is 'o', 
			in particular the Schur indices of $\chi $ at the infinite 
			places are 1. 
			Moreover in this case $q$ and hence $p$ is
			$\equiv 3\pmod{4} $ and ${\mathbb{Q}}_p$ does not contain 
			a $e_0$-th root of unity, so $\mathbb{Q}_p(\chi) $ is the 
			unramified extension of degree 2 of $\Q_p$. 
			As in (b) this field 
			is also the character field of all $p$-modular constituents
			of $\chi $, so the $p$-adic Schur index of $\chi $ is 1 
			by Theorem \ref{Benard}.
		For odd prime divisors $\ell $ of $p-1$, the character $\chi $
			remains irreducible modulo $\ell $ and hence again 
			all Schur indices are 1. 
		So the only prime where $\chi $ can have a 
		nontrivial 
		local Schur index is the unique prime of ${\mathbb{Q}}(\zeta_e^u )$ 
		that divides $2$. 
		As the sum of 
		the Hasse invariants is trivial, 
 all local Schur indices are 1.
	\end{itemize}

To compute the unitary discriminants of the characters in c) we use the 
strategy from \cite[Section 10]{Unitary} as described in Section \ref{OS}. 
Here we have 
 $q\equiv 3 \pmod{4}$ and $e_0\geq 4$.  The character field 
$L:=\mathbb{Q}(\mu )=\mathbb{Q}(\zeta_e^u)$ is a complex cyclotomic field  
of $2$-power order.
Let $K$ denote its maximal
real subfield. Let $\rho :A \to \GL_{p-1}(L) $ denote the representation
affording the character $\mu$. 

The group $A$ admits an automorphism $\alpha \in \Aut(A)$ with 
$$\alpha _{|Z} = \id _{|Z } , \alpha (\T_1 ) = \T_1 ^q  .$$ 
The restriction of $\alpha $ to the center of $A$ inverts all elements of 
$\langle \T_2 \rangle $ and hence $\alpha $ interchanges $\mu $ 
with its complex conjugate character. 
In particular, $K$ is the fixed field of the 
restriction of $\alpha $ to $L$. 
Put $\tilde{A} := A \rtimes  \langle \alpha \rangle $ to denote the semidirect 
product of $A$ with the group $\langle \alpha \rangle $ of order 2. 
Let $\widetilde{\rho } := \rho ^{\tilde{A}} $ be the induced representation and
$$ X:= \langle \rho (A) \rangle _L^{\alpha }  := \{ x\in  
\langle \rho (A) \rangle _L \mid \alpha (x)  = x \} $$ 
denote the fixed algebra of $\alpha $ in the enveloping algebra of $\rho (A)$.

Then by Theorem \ref{fixalg}  the class of $X$ is the class of the
enveloping algebra of $\widetilde{\rho}(\tilde{A})$ in the Brauer group of $K$. 
To determine $[X] \in \Br(K)$ we compute the local Schur indices of 
$\widetilde{\rho }$.
Now $\widetilde{\rho }$ is induced 
from the same linear character of the 
abelian normal subgroup $Z\times \langle \T _2 \rangle $ as $\rho $.

Recall that we are  in the case $q\equiv 3 \pmod{4} $ and put $c:=\frac{q-1}{2}$.
Then $c$ is odd and 
$\T_1^c \alpha $ acts as complex conjugation on ${\mathbb{Q}}(\exp(2\pi i/p),\zeta_e^u ) $ and satisfies 
	$$(\T_1^c\alpha )^2 =\T_1^c \T_1^{cq} \alpha ^2 = \T_1^{c(q+1)} = \T_0. $$
This allows to conclude that the real Schur index of $\widetilde{\rho} $ is 2 
	if and only if  $\widetilde{\rho}(\T_0) = -1 $,  so if and only if $u$ is odd.

	For the odd primes $\ell $ dividing $q-1$ the character of 
	$\widetilde{\rho }$ 
is in an $\ell $-block of defect $0$, so all $\ell $-local Schur indices are 1.
	As $2$ is totally ramified in the character field of $\widetilde{\rho }$,
the $2$-adic Schur index  can be read off from the
$p$-adic Schur indices. 
As above we can use \cite{Yamada} to pass to the group 
$A'\rtimes \langle \alpha \rangle $ whose Sylow $p$-subgroups are cyclic. 
As $\mathbb{Q}_p$ contains primitive $(p-1)$st roots of unity, 
 Theorem \ref{Benard} yields 
 that these $p$-adic Schur indices are all 1. 
So we get 
$$ {[}X{]}_K = \left\{ \begin{array}{ll} {[}K{]} & 
	\mbox{ if $u$ is even}  \\
{[}(-1,-1)_{K}{]} & \mbox{ if $u$ is odd.} \end{array} \right. $$

We are now in the position to apply Theorem \ref{fixalg}. 
From the above we obtain that 
$X$ is an orthogonal subalgebra
if and only if $u$ is even.
Then orthogonal determinant of the induced involution on $X$ can be obtained 
as the determinant of any skew symmetric unit in $X$ (see \cite[Proposition 2.2]{NebeOrtDet}).
Now $X$ 
 contains $\rho (Z) $ and the skew element 
$\rho (\z ) - \rho (\z ^{-1}) $ 
	has determinant $p^{(q-1)/(p-1)} \in q ({\mathbb{Q}}^{\times })^2 $.
By Theorem \ref{fixalg} we hence have $\disc(\mu ) = -q $ here. 

The Frobenius-Schur indicator of the character of $\widetilde{\rho }$ is -1, 
if and only if $u$ is odd, so here the restriction of the involution to $X$ is
symplectic. By Theorem \ref{fixalg} (c) the discriminant of $\mu $ is 
the $L$-discriminant of $[X]_K$. 
As $L=\mathbb{Q}(\zeta_e)$ contains a primitive fourth root of unity 
we have ${[}X{]}_K = (-1,-1)_K = (L,-1)_K$, so the 
$L$-discriminant of $[X]_K$ is $-1$.
\end{proof}

When computing unitary discriminants of  $\chi \in \Irr^{o}(G) $ for  $G\in \{ A_0, B , \SU_3(q) \}$ 
we will face the situation that the restriction of $\chi $ to $A$ contains 
a constituent from Lemma \ref{Adeg(q-1)} (b). 
Then Proposition \ref{ind+-} (a) shows that the contribution of this 
character to the unitary discriminant is trivial, if $q\equiv 1 \pmod{4}$. 
However, in the case where $q\equiv 3 \pmod{4}$ we need to compute the 
discriminant of $[{\mathcal Q}_p]$ over the character field of $\chi $. 
It turns out that in our situations $\Q(\chi )$ satisfies the assumption
of the next lemma, showing that the $\Q(\chi )$-discriminant of 
$[{\mathcal Q}_p]$ is $-p$. 

\begin{lemma} \label{quateq}
	Let $p$ be a prime, $p\equiv 3 \pmod{4}$. 
   Let $L$ be an abelian non-real number field with conductor dividing $p^a+1$ for 
	some odd integer $a$.
	Then $L$ splits ${\mathcal Q}_p$ 
	 and $\disc_L({\mathcal Q}_p) = -p$.
\end{lemma}

\begin{proof}
	Let $K$ denote the real subfield of $L$. We  show that 
	$(L,-p)_K \cong {\mathcal Q}_p \otimes K .$

	Let $\zeta $ be a primitive $(p^a+1)$th root of unity.
	Then $\Q_p(\zeta )$ is the unramified extension of degree $2a$ of 
	$\Q_p$. As the $a$th power of the Frobenius inverts $\zeta $ and
	hence is the complex conjugation on $\Q(\zeta )$,  all 
	$p$-adic completions of any non-real subfield of $\Q(\zeta )$ 
	have even degree over $\Q_p$. Moreover the $p$-adic completions of 
	the maximal real subfield 
	$\Q(\zeta + \zeta ^{-1} )$ have odd degree, $a$, over $\Q_p$. 

	In particular, $L$ splits ${\mathcal Q}_p$ as it splits this 
	algebra at the infinite places and at the unique ramified place, $p$. 
	Similarly we get that 
	${\mathcal Q}_p\otimes K$ is the quaternion algebra over $K$ 
	that is exactly  ramified at all infinite places of $K$ and 
 all the places of $K$ that divide $p$.  
	The same ramification behaviour holds for $(L,-p)_K$ 
	for the infinite primes and the ones dividing $p$. 

	It remains to show that no other primes ramify in 
	$(L,-p)_K$. 

	We first consider the dyadic primes $\lambda $ of $K$.
	If $p\equiv -1 \pmod{8}$ then $\sqrt{-p} \in \Q_2$, so 
	we can assume that $p$ and hence $p^a$ is congruent to $3$ mod $8$. 
	In particular, $\Q_2(\sqrt{-p})$ is the unramified quadratic extension
	of the $2$-adics. Moreover $p^a+1 = 4\cdot x$ for some odd number $x$ 
	and the completion $L_{(2)}$ of $L$ at $\lambda $ 
	is a subfield of $\Q_2(i) \Q_2 (\zeta _x)$, a field with inertia subfield
	$\Q_2(\zeta _x)$ and ramification degree $2$. 
	As $L_{(2)}/K_{(2)}$ is ramified of degree $2$, the completion 
	$K_{(2)}$ of $K$ at $\lambda $ is totally unramified of degree, say, 
	$f:=[K_{(2)} : \Q_2 ]$ over the $2$-adics. 
	If $f$ is even, then the unramified quadratic extension $\Q_2(\sqrt{-p})$ 
	is contained in $K_{(2)}$, so $(L,-p)_K$ is split at $\lambda $. 
	So $f$ is odd and the Galois group of $L_{(2)} / \Q _{(2)}$ is 
	abelian of order $2f$. In particular, 
	$$L_{(2)} = K_{(2)} \Q_2(\sqrt{a}) $$ is the compositum of 
	$K_{(2)}$ with a ramified quadratic extension $\Q_2(\sqrt{a}) $ of $\Q_2$ 
	and $$(L_{(2)},-p)_{K_{(2)}} = (a,-p)_{\Q_2} \otimes K_{(2)} .$$
	Now the conductor of $\Q_2(\sqrt{a})$ is not a multiple of $8$, so 
	$a$ is a unit in $\Z_2$ and hence a norm in the 
	unramified extension $\Q_2(\sqrt{-p})/\Q_2$. Therefore $(a,-p)_{\Q_2}$ is
	split and so is $(L_{(2)},-p)_{K_{(2)}}$.

	Now let $ \ell \neq p$ be an odd rational prime contained in some 
	prime $\lambda $ of $K$ that ramifies in $(L,-p)_K$. 
  The only prime ramifying in  $\Q (\sqrt{-p})/\Q $ is $p$.
	Therefore $\lambda $ ramifies in $L/K$, 
	and hence  $\ell $ divides the conductor of $L$, 
	so $\ell $ divides $p^a+1$, i.e. 
	 $(-p)^a \equiv 1 \pmod{\ell }$.  Recall that  $a$ is odd.
	 So $-p$ is a square mod $\ell $ and hence also in $\Q_{\ell }$. 
	Therefore primes dividing $\ell $ cannot ramify in $(L,-p)_K$. 
\end{proof}

From Lemma \ref{Adeg(q-1)} we now conclude Theorem \ref{A0deg(q-1)}: 

\begin{proof}  (of Theorem \ref{A0deg(q-1)}) 
	We have $A_0 = A \times \langle \T^{(q-1)e} \rangle $ and
	$\mu _{q-1}^{(u)}  = \mu \otimes \lambda $, 
	where $\mu = \Res_A(\mu _{q-1}^{(u)}) $ is unitary stable 
	and $\lambda $ 
 a linear character of the cyclic group $ \langle \T^{(q-1)e} \rangle $ of order $b$.
	If $u=0$ resp. $u=(q+1)/2$, then $\lambda = 1$. So case (a) and (b)
of Theorem \ref{A0deg(q-1)} follow from case (a) and (b) of Lemma \ref{Adeg(q-1)}.
In all other cases, the character field $L:=\Q(\mu_{q-1}^{(u)}) = \Q(\delta ^{(q-1)u}) $ 
	is a complex number field; in particular the Frobenius-Schur 
	indicator of $\mu _{q-1}^{(u)} $ is `o'. 
	If the Frobenius-Schur
        indicator of $\mu $ is `o' or $+$, then the unitary discriminant  of 
	$\mu_{q-1}^{(u)}$ is represented by any representative of $\disc(\mu )$, 
	for short 
	$$\disc (\mu_{q-1}^{(u)} )  = \disc(\mu \otimes \lambda ) = \disc(\mu )  .$$ 

	So it remains to consider the case where $u \equiv e/2 \pmod{e} $, i.e. 
	$[\mu] = [{\mathcal Q}_p ]$, but $u\neq (q+1)/2$. 
	For $q\equiv 3 \pmod{4} $ 
	the character field $L$ satisfies the assumption of Lemma \ref{quateq}.
	Using 
	Proposition \ref{ind+-} (a) we get that 
	$$\disc(\mu_{q-1}^{(u)} ) =  \left\{ \begin{array}{ll} 
		1  & q\equiv 1 \pmod{4} \\ 
	-q &  q\equiv 3 \pmod{4} .\end{array} \right. $$
\end{proof}

\subsection{The unitary discriminants of $B$} \label{Bunitary}

\begin{theorem}\label{Boreldet}
        The irreducible even degree characters of $B$ are the
	characters  $\vartheta_{q^2-q}^{(v)}$, $0\leq v\leq q$
	of degree $q(q-1)$ and the characters 
	$\vartheta_{(q^2-1)/d}^{(u,v)}$, $0 \leq u,v \leq d-1$. 
        \begin{itemize}
\item[a)] $\vartheta_{q^2-q}^{((q+1)/2)}$ is the character of a rational representation. Its orthogonal discriminant is $(-1)^{(q-1)/2}q$.
        \item[b)] $\vartheta_{q^2-q}^{(0)}$ 
is a rational character of Frobenius-Schur indicator $-$. Its
                Brauer element is $[\vartheta_{q^2-q}^{(0)} ] = [{\mathcal Q}_p]$.
                \item[c)]
                        For $v \not\in \{ 0, (q+1)/2 \} $ the Frobenius-Schur
			indicator of $\vartheta_{q^2-q}^{(v)}$  is 'o', the 
			character field is 
			$L= \Q(\vartheta_{q^2-q}^{(v)}) = \mathbb{Q}[\delta^{v(q-1)}] $ and 
			the unitary discriminant is 
			$$\disc(\vartheta_{q^2-q}^{(v)}) = (-1)^{(q-1)/2} q^{v-1} = 
			\left\{ \begin{array}{ll}
                -q & \mbox{ if }  (v,q) \equiv (0,3) \pmod{(2,4)} \\
                -1 & \mbox{ if } (v,q) \equiv (1,3) \pmod{(2,4)} \\
                 q & \mbox{ if }  (v,q) \equiv (1,1) \pmod{(2,4)} \\
                1 & \mbox{ if } (v,q) \equiv (0,1) \pmod{(2,4)}. 
                        \end{array} \right. $$
		\item[(d)]
	The characters $\vartheta_{(q^2-1)/d}^{(u,v)}$ have trivial unitary discriminant.
        \end{itemize}
\end{theorem}

\begin{proof}
{\bf (d)}
 The irreducible characters $\vartheta_{(q^2-1)/d}^{(u,v)}$ 
	of $B$ are trivial on $Z(U)$ and  induced 
from a non-trivial linear character  of $U/Z(U) $, an elementary abelian group
of order $q^2$. Let  $\T_0$  be the element of order $2$ in the torus $T$ 
	and consider the semidirect product $H := U\rtimes \langle \T_0 \rangle $, 
	a normal subgroup of $B$.
	Conjugation by $\T_0$ inverts the elements of $U/Z(U)$, so 
	the restriction $R:=\mathrm{Res}_H(\vartheta_{(q^2-1)/d}^{(u,v)})$ 
	of $\vartheta_{(q^2-1)/d}^{(u,v)}$  to $H$ 
	is the sum of orthogonal irreducible characters of 
	$H$ of degree 2. 
	The images of the corresponding degree 2 representations are 
	dihedral groups of order $2p$, so these constituents are orthogonal
	and of Schur index 1. 
 As the restriction of $R$ to $U$ is orthogonally stable, the orthogonal discriminant 
 of the unitary stable character real character $R$ is 1 by the 
	formula in \cite[Theorem 4.3 and Corollary 4.4]{Albanian}.
Moreover, no constituent of $R$ has a non-trivial Schur index,  so 
	we conclude that $[\vartheta_{(q^2-1)/d}^{(u,v)}]=1$ and 
	the unitary discriminant of $\vartheta_{(q^2-1)/d}^{(u,v)}$ 
is represented by 1.

	{\bf (a),(b),(c)} 
It remains to consider the characters $\vartheta_{q^2-q}^{(v)}$. 
These restrict non-trivially to $Z(U)$. 
 As $Z(U)$ is a normal subgroup of $B$, the trivial character of $Z(U)$
  does not occur in the restriction of $\vartheta_{q^2-q}^{(v)} $ to $Z(U)$, so
 the restriction of $\vartheta_{q^2-q}^{(v)} $ to the group $A_0$ 
	is a sum of the characters $\mu _{q-1}^{(u)}$  from Theorem \ref{A0deg(q-1)}.
From the character table of $B$ in \cite{Geck} we obtain the restriction of 
$\vartheta_{q^2-q}^{(v)}$ to the center of $A_0$ as 
	$$ \Res_{Z(A_0)}(\vartheta_{q^2-q}^{(v)}) = \sum_{u\neq v} (q-1) \lambda ^{(u)}  $$
	and hence 
	$$\Res_{A_0}(\vartheta_{q^2-q}^{(v)}) = \sum _{u\neq v} \mu_{q-1}^{(u)} .$$
	{\bf (a)} If $v=(q+1)/2$, then $\vartheta_{q^2-q}^{(v)}$ is a rational 
	orthogonal character. Its discriminant is given in \cite[Theorem 4.7]{SL3SU3}. 
	\\
	{\bf (b)} For $v=0$, $\vartheta_{q^2-q}^{(v)}$ is again rational, 
	but now its restriction to  $A_0$ contains the indicator $-$ character 
	$\mu_{q-1}^{((q+1)/2)} $ with multiplicity $1$ showing that the 
	Brauer element 
	$[\vartheta_{q^2-q}^{(0)}]= [\mu_{q-1}^{((q+1)/2)} ] = [{\mathcal Q}_p] $.
	\\
	{\bf (c)} In the remaining cases, the character field $L$ 
	of $\vartheta_{q^2-q}^{(v)}$  satisfies the assumptions of Lemma \ref{quateq}
	if $q\equiv 3 \pmod{4}$. 
	Moreover $R=\Res_{A_0}(\vartheta_{q^2-q}^{(v)})$ is unitary stable, so 
	we obtain the unitary discriminant of $\vartheta_{q^2-q}^{(v)} $
	from Theorem \ref{A0deg(q-1)}: 
	If $v$ is even, then $R$ contains 
		$(q+1)/2$ summands $\mu_{q-1}^{(u)} $ with odd $u$ and 
		$(q-1)/2$ summands $\mu_{q-1}^{(u)} $ with even $u$, so 
		$$\disc(\vartheta_{q^2-q}^{(v)} ) = \left\{ \begin{array}{ll} 
			-q & q\equiv 3 \pmod{4} \\
		1 & q\equiv 1 \pmod{4} . \end{array} \right. $$
	If $v$ is odd, then $R$ contains 
		$(q-1)/2$ summands $\mu_{q-1}^{(u)} $ with odd $u$ and 
		$(q+1)/2$ summands $\mu_{q-1}^{(u)} $ with even $u$,
	 so
		$$\disc(\vartheta_{q^2-q}^{(v)} ) = \left\{ \begin{array}{ll} 
			-1 & q\equiv 3 \pmod{4} \\
		q & q\equiv 1 \pmod{4} . \end{array} \right. $$
\end{proof}

For the characters 
	$\varphi^{(v)}$ and 
	$\tau ^{(u)}$ from Definition \ref{phiundtau} we find the following.

\begin{corollary} \label{detphiundtau}
	The unitary discriminant of $\tau ^{(u)} $ is trivial.
	Let $L$ be a complex number field so that 
	$L$ satisfies the assumption of Lemma \ref{quateq} in the 
	case where $q\equiv 3\pmod{4}$.
Then 
	$$ \disc_L(\varphi ^{(v)} ) = \left\{ \begin{array}{ll} 
		  q	 & \mbox{ if } q\equiv 1 \pmod{4} \\
	  1 & \mbox{ if } q\equiv 3 \pmod{4} \end{array} \right. 
		$$ 
\end{corollary}

\begin{proof} 
	All summands of $\tau ^{(u)}$ are of unitary discriminant 1 and hence
	$\disc(\tau ^{(u)}) = 1$. 
	\\
	As $L$ satisfies the assumptions of Lemma \ref{quateq} we see from Proposition \ref{ind+-} that 
	$$\disc_L(\vartheta _{q^2-q} ^{(0)}) = \left\{ \begin{array}{ll}
	1 & q \equiv 1 \pmod{4} \\ 
	-q & q \equiv 3 \pmod{4} \end{array} \right. $$
	and hence $\vartheta _{q^2-q} ^{(0)} $ contributes  in the same way 
	as all the other $\vartheta _{q^2-q} ^{(v)} $ with even $v$. 
	The character $\varphi ^{(u)}$ is a sum of $q+1$ irreducible characters 
	$\vartheta _{q^2-q}^{(v)} $ of which 
	$(q+1)/2$ have an odd upper index $v$ and 
	$(q+1)/2$ have an even upper index $v$. 
If $q \equiv 1 \pmod{4}$ then $(q+1)/2$ is odd and hence Theorem 
	\ref{Bunitary} yields $\disc_L(\varphi ^{(v)} ) = q $. 
	 For $q\equiv 3 \pmod{4}$ the number $(q+1)/2$ is even, so
	Theorem \ref{Bunitary} implies that 
	 $\disc_L(\varphi ^{(v)} ) = 1 $. 
\end{proof}

\subsection{The unitary discriminants of $\SU_3(q)$ for $q$ odd}

\begin{theorem} \label{mainSU3}
	Let  $q$ be a power of an odd prime $p$
	and put 
	$$f(u):= 
	\left\{ \begin{array}{ll} 
		-(\delta^{(q+1)u} - \delta^{-(q+1)u} )^2  & \mbox{ if } (q-1)/2 \mbox{ does not divide } u, \\ 
		-(\delta^{u} + \delta^{qu} )^2  & \mbox{ if } (q-1)/2 \mbox{ does divide } u. 
	\end{array} \right.
		$$ 
	The following table gives the unitary discriminant $\disc(\chi )$ 
	for the characters $\chi \in \Irr^o(\SU_3(q)) $. 

\begin{center}
\begin{tabular}{llll}
  \toprule
	$\chi$ & parameters &  $\disc(\chi )$ &  $\disc(\chi )$ \\
	& & $q \equiv 1 \pmod{4} $ & $q\equiv 3  \pmod{4}$ \\
  \midrule[0.06em]
                \addlinespace[0.3em]
	$\chi_{(q-1)(q^2-q+1)}^{(u,v,w)}$ & $\begin{array}{l} 
		1\leq u < v \leq (q+1)/d,\\ v<w\leq q+1 \\ u+v+w \equiv 0 \pmod{q+1}
	\end{array} $ &   $q$  &  $-q$  \\
 \midrule[0.06em]
                \addlinespace[0.3em]
	$ \chi_{q^3+1}^{( u )} 
	$ & 
	$\begin{array}{l}  1 \leq u < (q^2-1) \\ (q-1) \nmid u, (q+1) \nmid u \\(u) = (-uq) \end{array} $ &   
		$-q^{u+1}f(u) $ &
		$ q^{u+1}f(u)  $ \\
 \midrule[0.06em]
                \addlinespace[0.3em]
	$\chi_{(q+1)(q^2-1)}^{(u)}$ & $ \begin{array}{l} 1 \leq u < q^2-q+1 
	\\ (q^2-q+1)/d \nmid u \end{array} $ & 
	  $q$ &  $1$  \\
 \midrule[0.06em]
                \addlinespace[0.3em]
	$\chi_{(q+1)(q^2-1)/3}^{(u)}$ & $ 0 \leq u \leq 2$  & 
	  $q$  &  $1$ \\
 \midrule[0.06em]
                \addlinespace[0.3em]
	$\chi_{(q+1)(q^2-1)/3}^{(u)'}$ & $ 0 \leq u \leq 2$  & 
	  $q$  &  $1$ \\
\bottomrule[0.06em]
\end{tabular}
\end{center}
\end{theorem}

\begin{proof}
	The unitary discriminants of the characters $\chi $ of $\SU_3(q)$ are 
	obtained by restriction to the Borel subgroup $B$. 
We have $$\mathrm{Res}_B(\chi ) = \chi_T + \chi_U$$ where $\chi _U$ is unitary stable. 
	Note that for $q\equiv 3 \pmod{4} $ and all characters 
	$\chi \in \Irr^{o}(\SU_3(q)) $, the character field $L=\Q(\chi )$  
	satisfies the assumption of Lemma \ref{quateq}. 
	So we obtain $\disc(\chi _U)$ from Theorem \ref{Boreldet} 
using Proposition \ref{restoB}.

For all characters $\chi \in \Irr^o(\SU_3(q))$, except for 
the ones of degree $q^3+1$, the character $\chi_T$ is $0$, so 
it remains to consider the characters $\chi^{(u)}:=\chi_{q^3+1}^{(u)} $. 
Here $\chi_U^{(u)}(1) = q^3-1$ and 
$$\disc(\chi_U^{(u)} ) = (-1)^{q-1} q^{1+u}.$$
To handle the discriminant of the $U$-fixed space we use the 
	strategy and notation from Remark \ref{monirr}. 

	Let  $\Phi $ denote the Frobenius automorphism as in Section \ref{SU3}. 
Then $\chi^{(u)} \circ \Phi = \chi ^{(-u)}$ is the complex conjugate character
	and we are in the position to apply Theorem \ref{fixalg} 
	to compute the discriminant of  the 
	2-dimensional $J_U \Q \SU_3(q) J_U $-module
	$W^{(u)}:=\Q[\delta ^u] V^{(u)} J_U$. 
	As in Remark \ref{Wu}, a $\Q[\delta ^{u}]$-basis for $W^{(u)}$ is given by 
	$$(B , \sum _{\h \in U} B\w \h  ). $$
	We put $L:=\Q [\delta ^u + \delta ^{-qu} ] = \Q(\chi ^{(u)})$ 
	to denote the character field of $\chi^{(u)}$ and 
	let $K$ be its maximal real subfield. Then $K$ is also the 
	fixed field of $\Phi $ in $L$. 
	We denote the matrix representation of $J_U \Q\SU_3(q) J_U $ on $W^{(u)}$ by
	$\rho $ and  put $$R:=\langle \rho (J_U \g J_U ) \mid \g \in \SU_3(q) \rangle_L.$$
	Then $R$ is the $L$-algebra generated by the two matrices
	$$\rho(J_U \T J_U )=\mathrm{diag}(\delta ^u, \delta^{-qu }) \mbox{ and } 
	W:= \rho(J_U \w J_U )=\begin{pmatrix} 0 & 1 \\ 
	q^3 & 0 \end{pmatrix} .$$
	The form of $W $ can be obtained from explicit 
	computations in the
	Yokonuma algebra. 
	Alternatively, note that the second basis vector
	is the image of the first one and now the second row is obtained from the fact 
	that $W $ is self adjoint with respect to the 
	invariant form $\diag(q^3,1 )$. 

	The Frobenius automorphism $\Phi $ commutes with $J_U$, 
	fixes $\w $ and maps 
	$\T $ to $\T^q$. 
	The elements 
	$$x:=\T ^{(q+1)} - \T ^{-(q+1)}, \  y:= \T+\T^q-\T^{-1}-\T^{-q} $$ give rise to skew-symmetric elements  $X :=\rho(J_U x J_U)$ and $Y:=\rho(J_U y J_U)$ in
	the $\Phi $ fixed algebra $R^{\Phi }$.
If $u$ is not a multiple of $(q-1)/2$, then 
$$\det(X) = -(\delta ^{(q+1)u}-\delta ^{-(q+1)u}) ^2  =:
f(u) $$ 
is non-zero. 
	If $u$ is a multiple of $(q-1)/2 $, then  $2u/(q-1)$ is odd, 
	as $u$ is not a multiple of $(q-1)$.
	As $\delta ^{(q^2-1)/2} = -1$, we compute 
	$\delta ^{qu} = -\delta ^{-u} $ and hence 
	$$\det(Y) = -(\delta ^u + \delta^{qu} - \delta ^{-u} - \delta ^{-qu}) ^2 
	= -4 (\delta ^u + \delta ^{qu})^2 =: 
4 f(u) $$ 
is non-zero. 

To conclude that $\disc (\chi ^{(u)}_T) = -f(u)$  using Theorem \ref{fixalg}
it remains to show that $R^{\Phi} \cong K^{2\times 2}$. 
This is clear if $q$ is a square, since then the minimal polynomial 
of $W$ is reducible over $\Q $. 
So assume that $q$ is not a square.

If $u$ is not a multiple of $(q-1)/2$, then 
$R^{\Phi} $ is spanned as a $K$-algebra by $X_0:=2X-\mathrm{trace}(X) $ and $W$, so 
$$R^{\Phi }  \cong (\gamma ^2 , q^3 ) _K $$
where $\gamma =(\delta ^{u(q+1)}-\delta ^{-u(q+1)} ) \in \Q[\delta ]$. 
Then $K[X_0 ] \cong K[\gamma ] =   \Q [\delta ^{u(q+1)} ]$ is a maximal subfield of 
$R^{\Phi }$ and conjugation by $W$ induces the non-trivial Galois automorphism 
of $K[\gamma ]/K$. By Remark \ref{Ldiscdiv}, $R^{\Phi } \cong K^{2\times 2}$ 
if and only if no place $\wp $ of $K$ ramifies in $R^{\Phi}$. 
As $q^3>0$ all infinite places of $K$ are unramified in $R^{\Phi }$. 
By Remark \ref{Ldiscdiv2}, the finite places of $K$ that can possibly ramify 
in $R^{\Phi }$ are those dividing $q \disc(K[\gamma ] / K)$.
\\
The places $\wp $ dividing $q$ are split in $K[\gamma ]/K$, so they do not 
ramify in $R^{\Phi }$. Note that $K[\gamma ] = \Q[\delta^{u(q+1)}] $ is a
cyclotomic field and $K$ is its maximal real subfield. 
By \cite[Proposition 2.15]{Washington}, there are no finite ramified places in 
$K[\gamma ]/K$ unless $(q-1)/\gcd(2u,q-1)$ is a power of some prime $\ell $. 
In this case there is only one finite place $\wp$ of $K$ that is ramified
in $K[\gamma ]/K$. 
As the number of ramified places of $R^{\Phi }$ is even, 
also $R^{\Phi } \otimes K_{\wp } $ is split.

If $(q-1)/2$ divides $u$, then
$$R^{\Phi} = \langle Y/2 , W \rangle _{K} \cong (\gamma ^2 , q^3) $$ 
where $\gamma = \delta^u + \delta ^{qu}  = \delta^u-\delta^{-u}$.
So $K[Y] \cong K[\gamma ] \cong K[\delta ^u ] $ is a maximal subfield of 
$R^{\Phi }$ and conjugation by $W$ yields the non-trivial Galois automorphism 
of $K[\gamma ]/K$. 
Similarly as before, we conclude that $q$ is a norm in $K[\gamma ]/K$ and 
there is at mos one finite place of $K$ that is ramified in $K[\gamma ]/K$. 
As before this implies that 
$R^{\Phi }$ is split.  
\end{proof}

Note  that the proof above also shows that the 
Brauer elements $[\chi_{q^3+1}^{(u)} ] $ are trivial, a result that is also
obtained for all $\chi \in \Irr^o(\SU_3(q)) $ from Proposition \ref{restoB}. 
Note that the Schur indices for the irreducible characters of $\mathrm{PSU}_3(q)$ 
have been obtained by Gow \cite{GowSchur}, who also shows that all Schur indices 
of $\SU_3(q)$ divide $2$.

\bibliography{References}
\end{document}